\documentclass[reqno]{amsart}
\usepackage{amssymb}
\usepackage{amsmath}
\usepackage{amsfonts}
\usepackage{mathpazo}
\usepackage{graphicx}
\usepackage{subfig}
\usepackage[all]{xy}
\usepackage[lite,sorted-cites,compressed-cites]{amsrefs}
\numberwithin{equation}{section}
\input{tcilatex}

\newcommand{\R}{\mathbb{R}}
\newcommand{\C}{\mathbb{C}}

\newcommand{\sech}{\operatorname{sech}}
\DeclareMathOperator*{\Res}{Res}
\renewcommand{\b}{\textbf}
\renewcommand{\it}{\textit}

\begin{document}
\begin{abstract}
We discuss a new numerical schema for solving the initial value problem for the Korteweg-de Vries equation for large times. Our approach is based upon the Inverse Scattering Transform that reduces the problem to calculating the reflection coefficient of the corresponding Schr\"odinger equation. Using a step-like approximation of the initial profile and a fragmentation principle for the scattering data, we obtain an explicit recursion formula for computing the reflection coefficient, yielding a high resolution KdV solver. We also discuss some generalizations of this algorithm and how it might be improved by using Haar and other wavelets.
\end{abstract}

\title[KdV equation]{A Haar-type approximation and a new numerical schema for
the Korteweg-de Vries equation }
\author[UAF]{Jason Baggett \and Odile Bastille \and Alexei Rybkin}
\date{June, 2011}
\address{Department of Mathematics and Statistics \\
University of Alaska Fairbanks\\
PO Box 756660\\
Fairbanks, AK 99775
}
\email{jabaggett@alaska.edu}
\email{orbastille@alaska.edu}
\email{arybkin@alaska.edu}
\thanks{Based on research supported in part by the NSF under grants DMS 0707476 and DMS 1009673.}

\maketitle

\section{Introduction}
In this paper, we consider the well-known Korteweg-de Vries (KdV) equation
\begin{equation*}
	u_t - 6uu_x + u_{xxx} = 0
\end{equation*}
on the whole line with the initial condition $u(x,0) = V(x)$. The function $V(x)$ is assumed to be finite, nonpositive, and have compact support (i.e. zero outside of a finite interval). In particular, we will discuss a new algorithm for numerically approximating the KdV for large times $t$. For small times, there are many available algorithms to numerically integrate the KdV. Of particular interest is the operator splitting algorithm discussed in \cites{BV05, HKRT11}. These two algorithms can be coupled together to form a hybrid solver suitable for all times. 

The KdV equation is ``exactly solvable'' by relating it to the Schr\"odinger equation
\begin{equation*}
	-\phi_{xx} + V(x)\phi = \lambda \phi.
\end{equation*}
One obtains the so-called scattering data from the Schr\"odinger equation, and then the solution to the KdV can be obtained by performing the Inverse Scattering Transform (IST). In this sense, the IST linearizes the KdV (as well as some other nonlinear evolution PDEs). This provides us with an extremely powerful tool to analyze its solutions. Unfortunately, numerical algorithms based upon the IST are much less impressive and have not so far shown a noticeable improvement over conventional methods. The real power of the IST is in capturing the large time asymptotic behavior of solutions to the KdV equation (e.g. solitons) which is of particular interest in applications. As presented in \cite{AC91}, an asymptotic formula for the solution of the KdV can be obtained from the IST. Using this asymptotic formula, the large time solution of the KdV can be approximated from the scattering data alone without the full manchinery of the IST. This gives us a faster and more accurate method than other conventional methods for studying the large-time behavior of solutions to the KdV. Moreover, although our method is a PDE solver, it does not use any standard numerical PDE techniques; we need only calculate the scattering data using root finders and some linear algebra as described below.

The scattering data of the Schr\"odinger equation consists of the finitely many bound states $-\kappa_n^2$, the corresponding (left) norming constants $c_n^2$, and the (right) reflection coefficient $R$. The bound states are precisely the eigenvalues $\lambda$ of the Schr\"odinger equation that give square-integrable solutions $\phi$. The left and right reflection coefficients $L$ and $R$, respectively, and the transmission coefficient $T$ come from the asymptotic behavior of the left and right Jost solutions to the Schr\"odinger equation $\phi_l$ and $\phi_r$, respectively, where for $\lambda = k^2$
\begin{equation*}
	\phi_l(x,k) = \begin{cases} e^{ikx} + L(k)e^{-ikx} + o(1) & x \to -\infty	\\
										 T(k)e^{ikx} & x \to \infty.
					  \end{cases}
\end{equation*}
and
\begin{equation*}
	\phi_r(x,k) = \begin{cases} e^{-ikx} + R(k)e^{ikx} + o(1) & x \to \infty	\\
										 T(k)e^{-ikx} & x \to -\infty.
					  \end{cases}
\end{equation*}
If $\lambda = -\kappa^2$ is a bound state, then $\phi_l$ and $\phi_r$ are square-integrable. The corresponding left and right norming constants are defined by
\begin{align*}
	c_l = \left(\int_{-\infty}^\infty |\phi_l(x,i\kappa)T^{-1}(i\kappa)|^2 dx \right)^{-1} & & c_r =\left(\int_{-\infty}^\infty |\phi_r(x,i\kappa)T^{-1}(i\kappa)|^2 dx \right)^{-1}.
\end{align*}

With our assumption that $V$ has compact support, the bound states can be obtained as poles of $R$ in the upper-half complex plane, and the (left) norming constants can be retrieved from the residues at these poles. The poles of $R$ can be numerically approximated by using root finders. However, computing residues is numerically more difficult. We will instead consider a related function $B$ which is a rotation of the left reflection coefficient $L$. Then $B$ has the same poles as $R$, and its corresponding residues are equal to the residues of $R$ times a computable scaling factor. We give a new algorithm for computing the residues of $B$ as presented below.

We will approximate our potential using $N$ piecewise-constant functions. Then in each interval where the function is constant, the reflection and transmission coefficients $L_n$, $R_n$, and $T_n$ can be explicitly derived. Let
\begin{equation*}
	\Lambda = \begin{pmatrix}	1/T & -R/T	\\
															L/T & 1/\overline{T}	\\
							\end{pmatrix}.
\end{equation*}
Then the reflection and transmission coefficients $L$, $R$, and $T$ for the total potential can be derived from the principle of potential fragmentation (see, e.g. \cites{A92, AS02}), or layer stripping as it is known in the context of the Helmholtz equation \cite{SWG96}:
\begin{equation*}
	\Lambda = \Lambda_N ... \Lambda_2 \Lambda_1
\end{equation*}
for $k \in \R$ where bars denote complex conjugation and $\Lambda_n$ are the transition matrices
\begin{equation*}
		\Lambda_n = \begin{pmatrix}	1/T_n & -R_n/T_n	\\
															L_n/T_n & 1/\overline{T_n}	\\
							\end{pmatrix}.
\end{equation*}
This gives us a recursive formula for the left and right reflection coefficients and also for the function $B$. Using this recursive formula for $B$, we can derive a recursive matrix formula for the residues of $B$ at the poles in the upper-half complex plane. Consequently, all scattering data can be obtained, and the solution to the KdV can be numerically approximated for large times by the asymptotic formula given in \cite{AC91}.

In this paper, we also provide some numerical simulations. In particular, we give a comparison of computing bound states as poles of $R$ and $B$ and computing norming constants with our algorithm as opposed to other common algorithms. Although our algorithm is slower than standard methods for obtaining the scattering data, we demonstrate that it tends to be more accurate, especially for discontinuous initial profiles. We do not provide any error estimates; instead, the accuracy is verified on explicitly solvable examples. We also provide a comparison of the asymptotic solution to the KdV versus numerically integrated solutions.

Lastly, we also give some generalizations of our algorithm. For example, there is a natural generalization of our algorithm to higher order spline interpolants of $V(x)$. We also discuss possible improvements by using Haar and other wavelets. The Haar wavelets are piecewise constant functions that form an orthogonal system. Wavelets are closely related to Fourier series, and they exhibit many properties that are numerically desirable. Since we are approximating our potentials $V(x)$ using piecewise constant functions, one would believe that our algorithm can be modified to use Haar wavelets (and possibly more general wavelets).

\section{Notation}
We will denote the upper-half complex plane by $\C^+$. For a function $f : \C \to \C$, we will let $\overline{f}(z)$ denote complex conjugation and $\widetilde{f}(z) = f(-z)$. As is customary in analysis, we will let $L^2(\R)$ be the class of functions $f$ such that $\int_{\R} |f|^2 < \infty$. We will let $L_1^1(\R)$ denote the class of functions $f$ such that $\int_{\R} (1+|x|)|f(x)| < \infty$. Given functions $f,g \in L^2(\R)$, we define the $L^2$ inner product to be
\begin{equation*}
	\left\langle f,g \right\rangle_2 = \int_{\R} f\overline{g}.
\end{equation*}
and the $L^2$-norm $\| \cdot \|_2$ to be the norm with respect to this inner product, i.e.
\begin{equation*}
	\|f\|_2 = \left[\int_{\R} |f|^2 \right]^{1/2}.
\end{equation*}
For a set $A \subseteq \C$, $\chi_A$ will denote the characteristic function on $A$; i.e.
\begin{equation*}
	\chi_A(x) = \begin{cases} 1 & \text{if } x \in A	\\
									  0 & \text{otherwise}.
					\end{cases}
\end{equation*}

\section{Direct/Inverse Scattering Theory for the Schr\"{o}dinger Equation
on the Line}

Consider the KdV equation
\begin{equation*}
	u_t - 6uu_x + u_{xxx} = 0.
\end{equation*}
A particular stable solution of the KdV equation is given by 
\begin{equation*}
	u(x,t) = -2\kappa^2\sech^2(\kappa x - 4\kappa^3t + \gamma)
\end{equation*}
where $\kappa$ and $\gamma$ are real constants. Solutions of this form are called \it{solitons}. For more general initial profiles $u(x,0)$, in order to solve the KdV equation, we must first consider the Modified KdV equation (mKdV) 
\begin{equation*}
	v_t - 6v^2v_x + v_{xxx} = 0.
\end{equation*}
Miura (1967) discovered that one could obtain a solution to the KdV from a solution of the mKdV by the transformation $u = v_x + v^2$. Using Miura's transformation, we obtain
\begin{equation*}
	u_t - 6uu_x + u_{xxx} = \left(2v +\frac{\partial}{\partial x}\right)(v_t - 6v^2v_x + v_{xxx}).
\end{equation*}
Through translation and scaling, one can transform the equation
\begin{equation*}
	u_t - 6uu_x + u_{xxx} + \lambda u = 0
\end{equation*}
into the KdV
\begin{equation*}
	u_t - 6uu_x + u_{xxx} = 0.
\end{equation*}
Because of this, Miura's transformation takes the more general form $u - \lambda = v_x + v^2$.
If we assume that $v = \frac{\phi_x}{\phi}$, then Miura's transformation $u - \lambda = v_x + v^2$ becomes
\begin{equation*}
	-\phi_{xx}+u(x,t)\phi = \lambda \phi.
\end{equation*}
This equation is fundamental to Quantum Mechanics. It is known as the  one-dimensional, time-independent Schr\"odinger equation, where $\phi$ is the wave function, $u$ is the potential, and $\lambda$ is the energy. The problem of solving the KdV equation reduces to finding nontrivial solutions of the Schr\"odinger equation in $L^2(\R)$. However, not every $\lambda$ has such a solution. Hence, the Schr\"odinger equation is an eigenvalue problem. Moreover, the eigenvalues of the Schr\"odinger equation do not change over time. For this reason, we can replace $u(x,t)$ with our initial profile $u(x,0) = V(x)$, and solve
\begin{equation*}
	-\phi_{xx} + V(x)\phi = \lambda \phi
\end{equation*}
Suppose that $V \in L_1^1(\R)$. This ensures that there are finitely many soliton solutions.
We have then that $V(x) \to 0$ as $x \to \pm \infty$.  Hence, our solutions behave asymptotically like
\begin{equation*}
	-\phi_{xx} \sim \lambda \phi
\end{equation*}
Since $\phi$ is bounded, we must have that $\phi$ behaves asymptotically like a sinusoid ($\lambda > 0$) or decays like an exponential function ($\lambda < 0$).

Let $H = -\frac{d^2}{dx^2} + V(x)$. Then the Schr\"odinger equation becomes $H\phi = \lambda\phi$. The operator $H$ is called the \it{Schr\"odinger operator}. We have that $\lambda$ is an eigenvalue of the Schr\"odinger operator $H$ if $H -\lambda$ has no bounded inverse, and we say that $\lambda$ is in the \it{spectrum} of $H$. For each $\lambda > 0$, there is a nontrivial solution to $H\phi = \lambda\phi$. However, these eigenfunctions $\phi$ are not contained in $L^2(\R)$. We call the set of such $\lambda$ the \it{continuous spectrum} of $H$. The eigenvalues $\lambda < 0$ give square-integrable eigenfunctions $\phi$. However, there are only finitely many such $\lambda$. We call these $\lambda$ the \it{bound states}, while the set of bound states is called the \it{discrete spectrum}. The continuous spectrum gives rise to a component of the solution of the KdV which acts like a solution to the linear equation $u_t + u_{xxx} = 0$. This part of the solution is the dispersive portion of the wave. The discrete spectrum corresponds bijectively with the soliton solutions. This portion of the solution of the KdV is stable and does not decay over time. Thus, we are really only interested in knowing the discrete spectrum for large times.

\begin{figure}
	\centering
	\subfloat[Wave $e^{-ikx}$ radiating from $\infty$]{\includegraphics[height=2in]{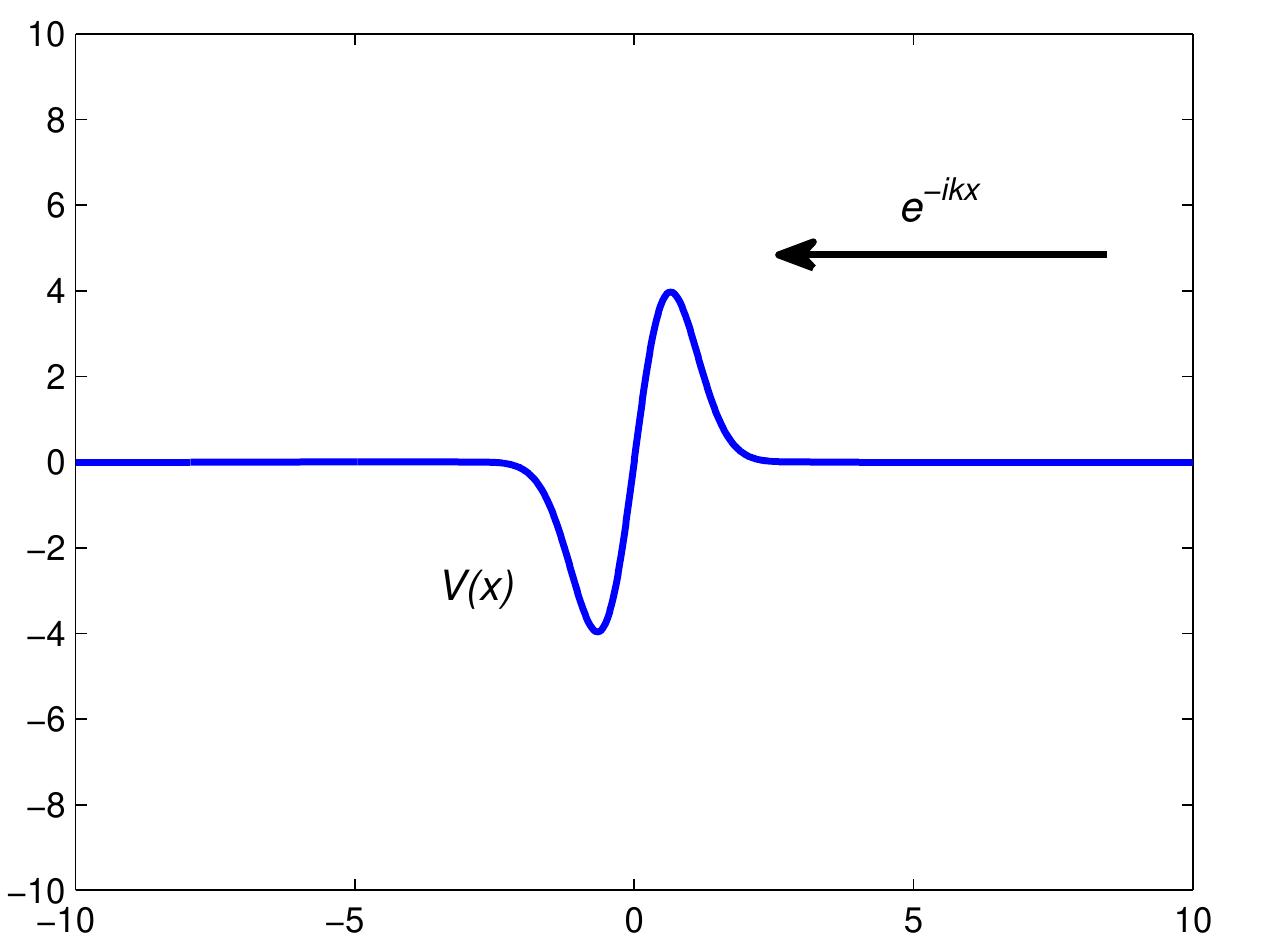}}
	\subfloat[$R(k)e^{ikx}$ is reflected]{\includegraphics[height=2in]{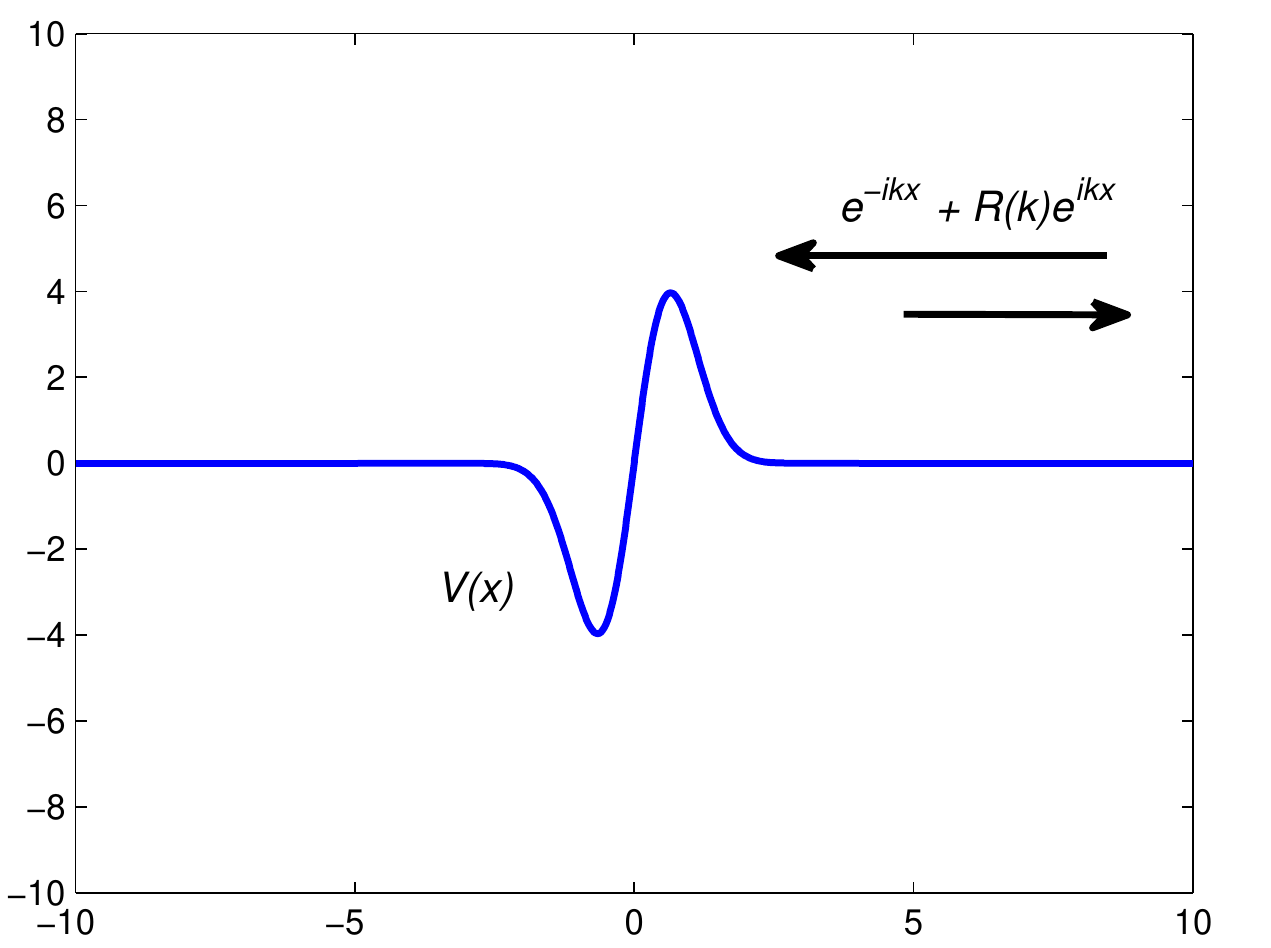}}	\\
	\subfloat[$T(k)e^{-ikx}$ is transmitted]{\includegraphics[height=2in]{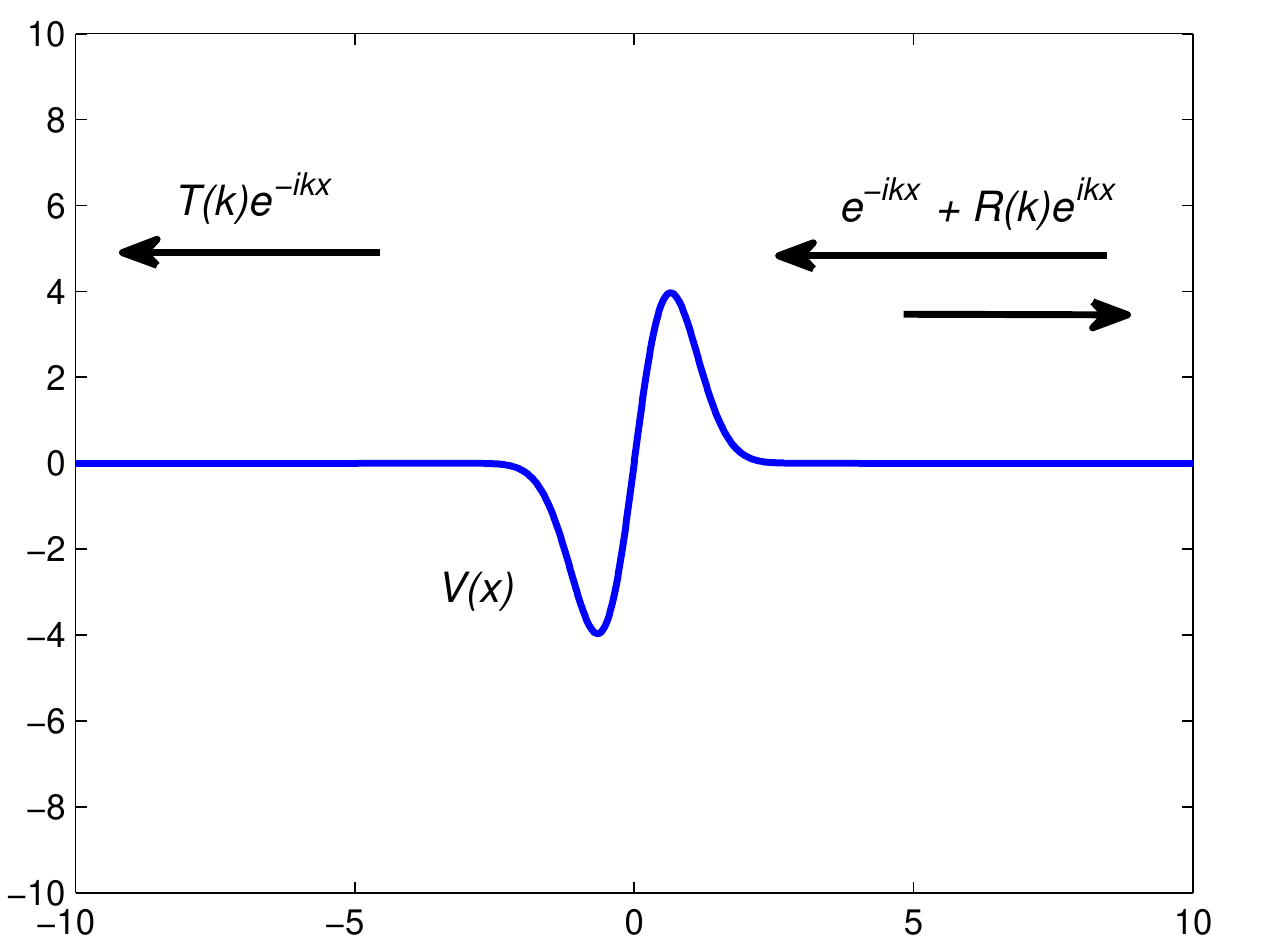}}
	\subfloat[Similarly, we can consider a wave coming from $-\infty$]{\includegraphics[height=2in]{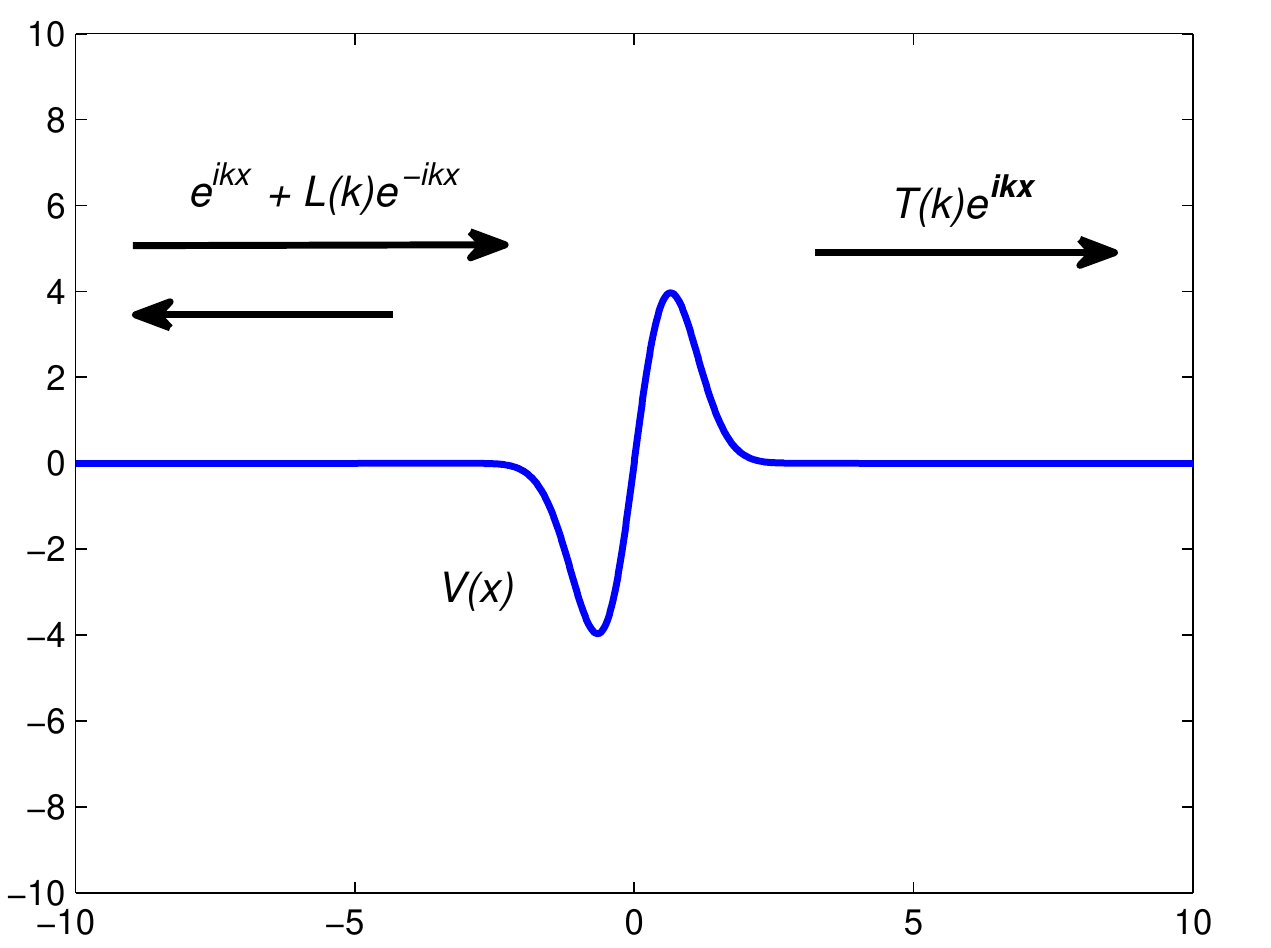}}
	\caption{The Jost solutions $\phi_l(x)$ and $\phi_r(x)$ as waves radiating from $\pm \infty$}
	\label{fig:wave}
\end{figure}

Suppose $\lambda = k^2 \in \R$. Among all solutions $\phi$ to the Schr\"odinger equation, pick the one satisfying 
\begin{equation*}
	\phi_r(x,k) = \begin{cases} e^{-ikx} + R(k)e^{ikx} + o(1) & x \to \infty	\\
										 T(k)e^{-ikx} & x \to -\infty.
					  \end{cases}
\end{equation*}
Such a solution $\phi_r$ is known as the \it{right Jost solution}. The function $T(k)$ is known as the \it{transmission coefficient}, and $R(k)$ is the \it{right reflection coefficient}. The reason for this terminology is that we can view $\phi_r$ as a wave $e^{-ikx}$ radiating from infinity, and $R(k)e^{ikx}$ is the portion of the wave that is reflected while $T(k)e^{-ikx}$ is the portion that is transmitted (see Figure \ref{fig:wave}). Similarly, we can consider the \it{left Jost solution}
\begin{equation*}
	\phi_l(x,k) = \begin{cases} e^{ikx} + L(k)e^{-ikx} + o(1) & x \to -\infty	\\
										 T(k)e^{ikx} & x \to \infty.
					  \end{cases}
\end{equation*}
where $T(k)$ is the same transmission coefficient and $L(k)$ is the \it{left reflection coefficient}.

If $k = i\kappa$ is a bound state, then $\phi_l(x,i\kappa), \phi_r(x,i\kappa) \in L^2(\R)$. We define the \it{left} and \it{right norming constants} at $k = i\kappa$ to be
\begin{align*}
	c_l = \|\phi_l(x,i\kappa)T^{-1}(i\kappa)\|_2^{-1} & & c_r = \|\phi_r(x,i\kappa)T^{-1}(i\kappa)\|_2^{-1}
\end{align*}

\section{The Classical Inverse Scattering Transform}
%Put here the relevant stuff up to page 20 of Jason's talk. 
Since $V(x) \in L_1^1(\R)$, we have that there are finitely many bound states $\lambda = k^2$ where $k = i\kappa$. Let $K$ denote the number of bound states, and let
\begin{equation*}
	\kappa_1 > \kappa_2 > ... > \kappa_K > 0.
\end{equation*}
Let $c_n$ denote the left norming constant at $k = i\kappa_n$.

Once we know the scattering data for the Schr\"odinger operator, we can use the Inverse Scattering Transform (IST) to obtain the soliton solutions of the KdV equation.
\centerline{
  \xymatrix{u(x,0)\ar[rrr]^{\text{direct scattering}} &&&S(0)\ar[d]^{\text{time evolution}}\\
    u(x,t) &&&S(t)\ar[lll]^{\text{inverse scattering}}\\
  } }
The Direct Scattering Transform is we map the initial potential $u(x,0)$ into the scattering data
\begin{equation*}
	S(0) = \{\{-\kappa_n^2,c_n\}_{n=1}^{K},R(k), k \in \R\}.
\end{equation*}
Next, we evolve the scattering data over time in a simple fashion:
\begin{itemize}
\item $\kappa_n(t) = \kappa_n$
\item $c_n(t) = c_ne^{4\kappa_n^3t}$
\item $R(k,t) = R(k)e^{8ik^3t}$
\end{itemize}
where $\kappa_n = \kappa_n(0)$, $c_n = c_n(0)$, and $R(k) = R(k,0)$. Then the scattering data becomes
\begin{equation*}
	S(t) = \{\{-\kappa_n(t)^2,c_n(t)\}_{n=1}^{K},R(k,t), k \in \R\}.
\end{equation*}

We can reclaim the solution to the KdV using Inverse Scattering as follows:
\begin{itemize}
\item
Form the Gelfand-Levitan-Marchenko (GLM) kernel:
\begin{equation*}
	F(x,t) = \sum_{n=1}^{N}c^2_{n}(t)e^{-\kappa _{n}x}+\frac{1}{2\pi }\int_{-\infty}^{\infty}e^{ikx}R(k,t)dk.
\end{equation*}
\item 
Solve the Gelfand-Levitan-Marchenko equation for $K(x,y,t)$, $y \geq x$: 
\begin{equation*}
K(x,y,t) + F(x+y,t) +\int_{x}^{\infty}F(s+y,t) K(x,s,t)ds = 0.
\end{equation*}
\item 
The solution to the KdV equation is 
\begin{equation*}
	u(x,t) = -2\frac{\partial}{\partial x}K(x,x,t).
\end{equation*}
\end{itemize}

Luckily, for large times $t$ we can simplify the GLM kernel. We have that 
\begin{equation*}
	\int_{-\infty}^{\infty}e^{ikx}R(k,t)dk \to 0 \text{ as } t \to \infty.
\end{equation*}
Thus, for large times we can approximate the GLM kernel by
\begin{equation*}
	F(x,t) \approx \sum_{n=1}^{N}c^2_{n}(t)e^{-\kappa _{n}x}.
\end{equation*}
Let
\begin{equation*}
	C(x,0) = \begin{pmatrix} c_{11}(x) & c_{12}(x) & \ldots & c_{1,N}(x) \\
											 c_{21}(x) & c_{22}(x) & \ldots & c_{2,N}(x) \\
											 \vdots & & \ddots & 		\\
											 c_{N,1}(x) & c_{N,2}(x) & \ldots & c_{N,N}(x)
			 	 \end{pmatrix}
\end{equation*}
where
\begin{equation*}
	c_{mn}(x) = \frac{c_{m}c_{n}}{\kappa_{m}+\kappa_{n}}e^{-(\kappa_{m}+\kappa _{n}) x}.
\end{equation*}
The matrix $C$ evolves in time by
\begin{equation*}
	c_{mn}(x,t) = c_{mn}(x) e^{4(\kappa_{m}^{3}+\kappa _{n}^{3}) t}
\end{equation*}
Then for large times, our solution to the KdV is \cites{AC91}
\begin{equation*}
	u(x,t) \approx -2\frac{\partial ^{2}}{\partial x^{2}}\ln [\det(I+C(x,t))]
\end{equation*}
From this, one obtains the asymptotic formula \cite{AC91}
\begin{equation*}
	u(x,t) \sim -2 \sum_{n=1}^N \kappa_n^2 \sech^2(\kappa_nx - 4\kappa_n^3t + \ln\sqrt{\gamma_n})
\end{equation*}
where
\begin{equation*}
	\gamma_n = \frac{2\kappa_n}{c_n^2}\prod_{m=1}^{n-1}\left(\frac{\kappa_n+\kappa_m}{\kappa_n-\kappa_m}\right)^2.
\end{equation*}
Notice that the large time solution $u(x,t)$ of the KdV behaves like a finite sum of single solitons. Moreover, we no longer need to do the full IST to solve the KdV for large times. We need only find the bound states $-\kappa_n^2$ and norming constants $c_n$.

If $R$ and $T$ can be analytically continued, then the poles of $R$ and $T$ in $\C^+$ are precisely $i\kappa_n$. That is, all of the poles of $R$ and $T$ in $\C^+$ lie on the imaginary axis and correspond with the bound states. Better yet, these poles are actually simple \cites{A94, MD05}. 
Furthermore, if we assume that $V$ has compact support, then \cite{AK01}
\begin{equation}
	\Res_{k=i\kappa_n}R(k) = ic_n^2.	\label{ResR}
\end{equation}
Consequently, if our potential $V$ has compact support and $R$ is analytically continued into $\C^+$, then the bound states and norming constants can be obtained from knowledge of $R(k)$ for $k \in \C^+$. In this case, we can approximate the solution of the KdV for large times from only knowledge of $R(k)$ for $k \in \C^+$.

\section{The scattering Quantities for a Block (Well) Potential}

\begin{figure}
	\centering
	\includegraphics[width=.9\textwidth]{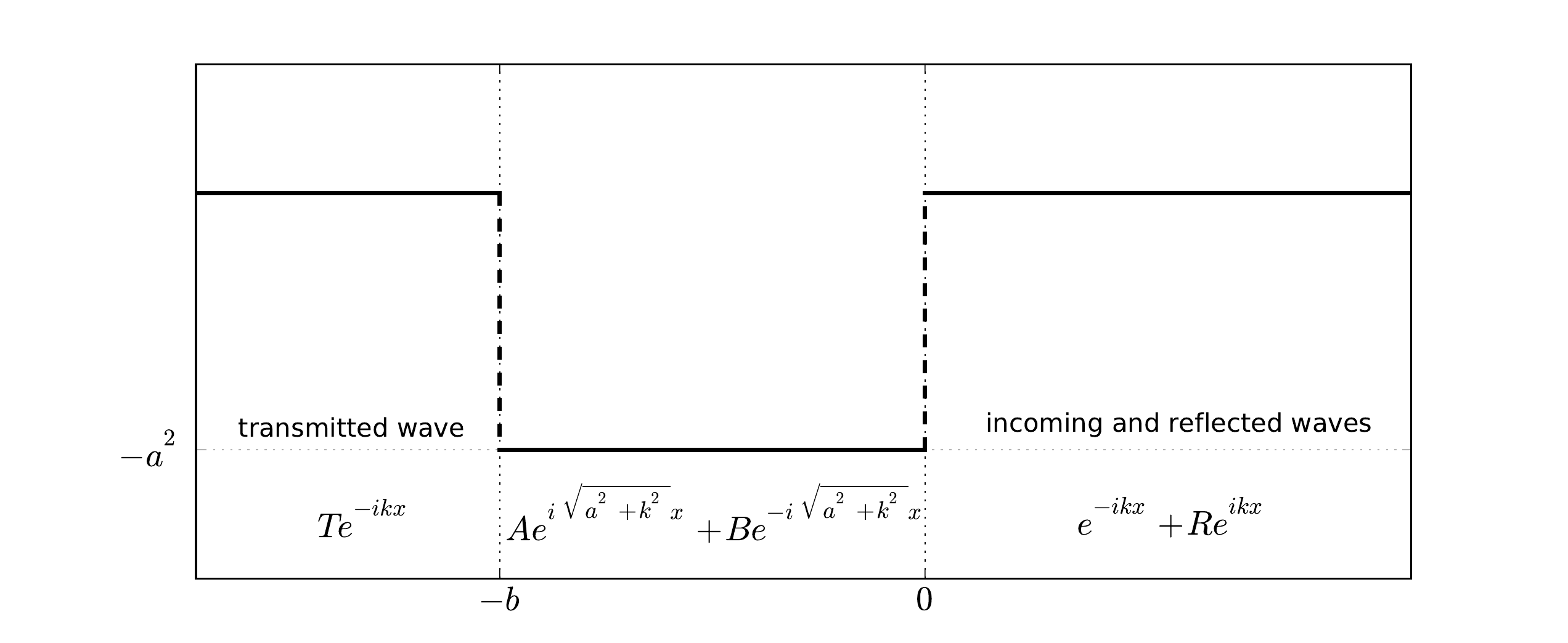}
	\caption{The setup for a single block potential}
	\label{fig:one_block}
\end{figure}

Consider the case when our potential $V$ is a single nonpositive well which is $-a^2$ on the interval $[-b,0]$ and 0 elsewhere, i.e. $V(x) = -a^2\chi_{[-b,0]}(x)$ (see Figure \ref{fig:one_block}).
In this case, we can obtain an exact solution to the Schr\"odinger equation. Moreover, using the continuity of the solution $\phi$ and its derivative $\phi_x$, we can set up a system of equations and solve for $R$ and $T$. Doing this, we obtain
\begin{equation}
	R(k)= \omega^2 \frac{1-\xi}{\xi-\omega^4} \quad  , \quad L(k) = \omega^2 \frac{1-\xi}{\xi-\omega^4}e^{-iab(\omega-1/\omega)} \quad , \quad T(k)=\frac{1-\omega^4}{\xi-\omega^4}e^{i\frac{ab}{\omega}}	\label{oneblock}
\end{equation}
where
\begin{equation*}
	\omega=\frac{k}{a}+\sqrt{\left(\frac{k}{a}\right)^2+1} \qquad, \qquad \xi=e^{i(\omega+1/\omega)ab}. 
\end{equation*}
These formulas for $R$, $L$, and $T$ are actually meromorphic in $\C$ if we choose the branch cut along the imaginary axis between $-ia$ and $ia$. Using these formulas, $R$, $L$, and $T$ can be analytically continued in $\C^+$. The only difficulty lies in considering the branch cut. However, we have that $\omega(-\overline{k}) = -\overline{\omega(k)}$ and $\xi(-\overline{k}) = \overline{\xi(k)}$. It follows that $R(-\overline{k}) = \overline{R(k)}$ and $T(-\overline{k}) = \overline{T(k)}$. For $k \in i\R$, we have that $R(k) = R(-\overline{k}) = \overline{R(k)}$ and $T(k) = T(-\overline{k}) = \overline{T(k)}$, so $R$ and $T$ are real-valued for $k \in i\R$. For $k = +0+i\kappa$, we have that $-\overline{k} = -0+i\kappa$. Therefore, $R(-0 + i\kappa) = \overline{R(+0 + i\kappa)}$ and since $R$ is real-valued on $i\R$, $\overline{R(+0+i\kappa)} = R(+0+i\kappa)$. Hence, $R(-0 + i\kappa) = R(+0 + i\kappa)$, so $R$ is continuous along the branch cut between $-ia$ and $ia$. It follows that $R$ is meromorphic in $\C$. Similarly, $T$ is meromorphic in $\C$ as well.

Consider the poles $i\kappa_n$ and residues $ic_n^2$ of $R$. The poles of $R$ and $T$ satisfy $\xi = \omega^4$. If we let $y_n = \frac{\kappa_n}{a}$, then $\kappa_n$ and $c_n$ can be explicitly computed by the following formulas:
\begin{equation}
	\frac{ab}{\pi}\sqrt{1-\left(\frac{\kappa_n}{a}\right)^2}-\frac{2}{\pi}\arctan \frac{\kappa_n}{a\sqrt{1-\left( \frac{\kappa_n}{a} \right)^2}}= n - 1 
	\label{bound_oneblock}
\end{equation}
and
\begin{equation}
	c_n^2=\frac{2\kappa_n\left(1-\left(\frac{\kappa_n}{a}\right)^2\right)}{2+b\kappa_n}	\label{norm_oneblock}
\end{equation}
for $n = 1,...,\left\lceil \frac{ab}{\pi}\right\rceil$.

\section{The Potential Fragmentation and the Scattering Quantities for Potentials Composed of Blocks}	\label{sec:multiblocks}
We define the \it{scattering matrix} to be
\begin{equation*}
	S = \begin{pmatrix} T & R \\ L & T \end{pmatrix}
\end{equation*}
The matrix $S$ is \it{unitary}, i.e. $S^{-1} = S^*$ where $S^*$ is the conjugate transpose of $S$ \cite{A94}. This gives us a few identities, namely \cites{AKM96, MD05}
\begin{equation}
	L\widetilde{T} + T\widetilde{R} = 0 \label{LTTR}
\end{equation}
for $k \in \R$. 
If we were to shift our potential to the right by $p$, then the scattering matrix would change as follows \cite{A92}:
\begin{align}
	L(k - p) &= L(k) e^{2ikp}	 \label{shiftL}	\\
	T(k - p) &= T(k)					\label{shiftT}	\\
	R(k - p) & = R(k) e^{-2ikp}	\label{shiftR}
\end{align}

Now suppose that our potential $V$ consists of $N$ nonpositive blocks. Let $R_n, L_n, T_n$ be the reflection and transmission coefficients on the $n$-th block: $V_n(x) = -a_n^2$ on $[-b_n, -b_{n-1}]$ where $b_0 = 0$. Let $R_n^0, L_n^0, T_n^0$ be the reflection and transmission coefficients on the $n$-th block shifted to the origin: $V_n(x) = -a_n^2$ on $[-(b_n-b_{n-1}),0]$. Let $R_{1,2,...,n}$, $L_{1,2,...,n}$, $T_{1,2,...,n}$ be the reflection and transmission coefficients on the first $n$ blocks. $R,L,T$ without subscripts or superscripts will denote the reflection and transmission coefficients for the overall potential.

Let
\begin{equation}
	\Lambda = \begin{pmatrix}	1/T & -R/T	\\
															L/T & 1/\widetilde{T}	\\
							\end{pmatrix}.	\label{transition_matrix}	
\end{equation}
The \it{fragmentation principle}, or \it{layer stripping principle} as it is also known, \cites{A92, AS02, AKM96, SWG96} says that for $k \in \R$
\begin{equation}
	\Lambda = \Lambda_N ... \Lambda_2 \Lambda_1	\label{potentialfragmentation}
\end{equation}
where $\Lambda_n$ are the \it{transition matrices}
\begin{equation*}
		\Lambda_n = \begin{pmatrix}	1/T_n & -R_n/T_n	\\
															L_n/T_n & 1/\widetilde{T_n}	\\
							\end{pmatrix}
\end{equation*}
(Note that blocks with $a_n = 0$ may be simply ignored since this implies $\Lambda_n$ is the identity matrix).
We can use potential fragmentation to come up with some recursive formulas. Using \ref{LTTR}, \ref{shiftL}, \ref{shiftT}, \ref{shiftR} and \ref{potentialfragmentation}, we obtain that
\begin{equation}
R_{1,...,n+1} = -\frac{L_{1,...,n}}{\widetilde{R_{1,...,n}}} \frac{R_{n+1}^0e^{2ikb_n}-\widetilde{L_{1,...,n}}}{1-R_{n+1}^0L_{1,...,n}e^{2ikb_n}}. \label{fragR}
\end{equation}
A similar expression may be obtained for the left reflection coefficient:
\begin{equation}
L_{1,...,n+1} = -\frac{R_{n+1}^0}{\widetilde{R_{n+1}^0}} \frac{L_{1,...,n}e^{2ikb_n}-\widetilde{R_{n+1}^0}}{1-R_{n+1}^0L_{1,...,n}e^{2ikb_n}} e^{-2ikb_{n+1}}. \label{fragL}
\end{equation}
We have that $L_{1,...,n} = -\frac{T_{1,..,n}}{\widetilde{T_{1,...,n}}}\widetilde{R_{1,...,n}}$ for $k \in \R$. Thus, $\left|L_{1,...,n}\right| = \left|R_{1,...,n}
\right|$ for $k \in \R$.
Therefore, $L_{1,...,n} = R_{1,...,n}e^{-2ik\beta_n}$ for some $\beta_n: \R \to \R$.
Equations \eqref{fragR} and \eqref{fragL} then give us that
\begin{equation}
R_{1,...,n+1} = -\frac{R_{1,...,n}}{\widetilde{R_{1,...,n}}} \frac{R_{n+1}^0e^{2ik(b_n-\beta_n)}-\widetilde{R_{1,...,n}}}{1-R_{n+1}^0R_{1,...,n}e^{2ik(b_n-\beta_n)}}. \label{fragR2}
\end{equation}
where
$\beta_1 = b_1$ and
\begin{equation}
e^{-2ik\beta_{n+1}} = \frac{R_{n+1}^0}{\widetilde{R_{n+1}^0}}\frac{\widetilde{R_{1,...,n}}}{R_{1,...,n}} 
\frac{R_{1,...,n}e^{2ik(b_n-\beta_n)}-\widetilde{R_{n+1}^0}}{R_{n+1}^0e^{2ik(b_n-\beta_n)}-\widetilde{R_{1,...,n}}}e^{-2ikb_{n+1}}. \label{expbetaN}
\end{equation}
Define $A_n = \frac{L_{1,...,n}}{R_{1,...,n}}e^{2ikb_n}$. Then $A_n = e^{2ik(b_n-\beta_n)}$ for $k \in \R$.
Equations \eqref{fragR2} and \eqref{expbetaN} give us that
\begin{equation}
R_{1,...,n+1} = -\frac{R_{1,...,n}}{\widetilde{R_{1,...,n}}} \frac{A_nR_{n+1}^0-\widetilde{R_{1,...,n}}}{1-A_nR_{n+1}^0R_{1,...,n}} \label{recR}
\end{equation}
and 
\begin{equation}
A_{n+1} = \frac{R_{n+1}^0}{\widetilde{R_{n+1}^0}}\frac{\widetilde{R_{1,...,n}}}{R_{1,...,n}}  
\frac{A_nR_{1,...,n}-\widetilde{R_{n+1}^0}}{A_nR_{n+1}^0-\widetilde{R_{1,...,n}}}	\label{recA}
\end{equation}
where $A_1 = 1$.
Let us next define $B_n = A_nR_{1,...,n} = L_{1,...,n}e^{2ikb_n}$. Then we get the following recursive formula:
\begin{equation}
	B_{n+1} = -\frac{R_{n+1}^0}{\widetilde{R_{n+1}^0}}\frac{B_n-\widetilde{R_{n+1}^0}}{1-R_{n+1}^0B_n} \label{recB}
\end{equation}
where $B_1 = R_1$.
Notice that $B_{n+1}$ is a M\"obius transform of $B_n$, and that the recursive formula for $B_n$ is much simpler than the recursive formula for $R_{1,...,n}$. Moreover, this formula only depends on $B_n$, $R_{n+1}^0$, and $\widetilde{R_{n+1}^0}$. From \ref{oneblock},
\begin{equation}
R_n^0 = \frac{\omega_n^2(1-\xi_n)}{\xi_n - \omega_n^4}. \label{Rn0}
\end{equation}
where $h_n = b_n - b_{n-1}$ is the width for the $n$-th block, $\omega_n(k) = \frac{k}{a_n} + \sqrt{(\frac{k}{a_n})^2 + 1}$, and $\xi_n(k) = e^{-ia_nh_n(\omega_n(k)+1/\omega_n(k))}$. For $k \in \R$, we have that $\overline{R_n^0} = \widetilde{R_n^0}$. By taking the complex conjugate of \eqref{Rn0}, we obtain that for $k \in \R$,
\begin{equation}
	\widetilde{R_n^0} = \frac{\omega_n^2(1-\xi_n)}{\xi_n \omega_n^4 - 1}. \label{CRn0}
\end{equation}
Since $R_n^0$ is meromorphic in $\C$, $\widetilde{R_n^0}$ is meromorphic in $\C$ as well (in particular, both are meromorphic in $\C^+$ where the poles of interest lie). Since the formula in \eqref{CRn0} is meromorphic in $\C$, it follows that \eqref{CRn0} holds for all $k \in \C$. Continuing inductively using equations \eqref{recR}, \eqref{recA}, \eqref{recB}, it follows that $R_{1,...,n}, A_n,$ and $B_n$ can be continued to meromorphic functions in $\C$ (and in particular, $\C^+$) for all $1 \leq n \leq N$.

Since $B_n = A_nR_{1,...,n} = L_{1,...,n}e^{2ikb_n}$, we have that $B_n$ has the same poles $k = i\kappa_m$ in $\C^+$ as $L_{1,..,n}$. Consequently, $B_n$ and $R_{1,...,n}$ have the same poles in $\C^+$. Since the poles $k = i\kappa_m$ in $\C^+$ of $L_{1,...,n}$ and $R_{1,...,n}$ are simple, we have that $A_n = \frac{L_{1,...,n}}{R_{1,...,n}}e^{2ikb_n}$ is analytic in $\C^+$ and nonzero at all $k = i\kappa_m$. It follows from equation \eqref{ResR} then that
\begin{align}
	\Res_{k = i\kappa_m}B_n &= A_n(i\kappa_m)\Res_{k = i\kappa_m}R_{1,...,n}  \\
													&= ic_m^2 A_n(i\kappa_m)	\label{resR2}
\end{align}
The value of $A_n(i\kappa_m)$ can be determined via the recursive formula \ref{recA}. If we can determine an algorithm for determining the residues of $B_N$, then this would effectively give us an algorithm for calculating the (left) norming constants.

Now suppose $B_n$ has the form
\begin{equation}
B_n = -\frac{R_n^0}{\widetilde{R_n^0}} \frac{p_n}{q_n}. \label{recB2}
\end{equation}
Applying \eqref{recB}, we get a linear system of recurrence relations for $p_n$ and $q_n$:
\begin{align}
\begin{array}{lcrcr}
p_{n+1} &=& -R_{n}^0p_n& - &\widetilde{R_{n+1}^0}\widetilde{R_n^0}q_n \label{pqrec}
\\q_{n+1} &=& R_{n+1}^0R_n^0p_n& + &\widetilde{R_n^0}q_n
\end{array}
\end{align}
or in matrix form
\begin{equation}
\begin{pmatrix} p_{n+1} \\ q_{n+1} \end{pmatrix} = M_n \begin{pmatrix} p_n \\ q_n \end{pmatrix} = M_n...M_2M_1\begin{pmatrix}p_1 \\ q_1 \end{pmatrix} \label{mateq}
\end{equation}
where
\begin{equation*}
M_i = \begin{pmatrix} -R_i^0 & -\widetilde{R_{i+1}^0}\widetilde{R_i^0}
										\\ R_{i+1}^0R_i^0 & \widetilde{R_i^0} \end{pmatrix}.
\end{equation*}

Let $N$ denote the number of blocks. If $q_N(k) = 0$ but $k$ is not a pole of $B_N$, then $p_N = 0$ as well. From \eqref{mateq}, this means that $\det(M_n) = 0$ for some $1 \leq n \leq N-1$ or $\begin{pmatrix} p_1 \\ q_1\end{pmatrix} = \b{0}$.
Since $B_1 = R_1$, from \eqref{recB2} we have that $\frac{p_1}{q_1} = -\widetilde{R_1}$. Our choice of $p_1$ and $q_1$ is arbitrary, so long is this ratio is preserved, since our resulting solution of $B_{n+1}$ is independent of our choice for $p_1$ and $q_1$. Some choices for our initial vector may be preferable for numerical computations, but for our purposes we will choose $\begin{pmatrix} p_1 \\ q_1\end{pmatrix} = \begin{pmatrix} -\widetilde{R_1} \\ 1 \end{pmatrix}$, because it is nonzero for all $k$. Thus, if $q_N = 0$ but $k$ is not a pole of $B_N$, then $\det(M_{N-1}...M_2M_1) = 0$. Equivalently, if $q_N = 0$ and $\det(M_{N-1}...M_2M_1) \neq 0$, then $k$ is a pole of $R_{1,...,N}$.

We claim that $\det(M_{N-1}...M_2M_1)(k) = 0$ for some $k \in \C^+$ if and only if $k = ia_N$ or for some $1 \leq n \leq N-1$ and some $0 \leq m \leq \left\lfloor \frac{a_n h}{\pi} \right\rfloor$,
\begin{equation}
	k = i\sqrt{a_n^2 - \left( \frac{\pi m}{h_n} \right)^2}. \label{detM}
\end{equation}.
We have that $\det(M_{N-1}...M_2M_1) = 0$ if and only if $\det(M_n) = 0$ for some $1 \leq n \leq N-1$. Moreover, 
\begin{equation*}
	\det(M_n) = R_n^0\widetilde{R_n^0}(R_{n+1}^0\widetilde{R_{n+1}^0}-1).
\end{equation*}
Thus, $\det(M_n) = 0$ if and only if $R_n^0 = 0$ (equivalently $\widetilde{R_n^0} = 0$) or $R_{n+1}^0\widetilde{R_{n+1}^0} = 1$. The second case occurs when 
\begin{equation*}
	\omega_{n+1}^4(1-\xi_{n+1})^2 = (\xi_{n+1}-\omega_{n+1}^4)(\xi_{n+1}\omega_{n+1}^4 -1).
\end{equation*}
After some algebra and noting that $\xi_{n+1} = e^{-ia_{n+1}h_{n+1}(\omega_{n+1}+1/\omega_{n+1})} \neq 0$, this simplifies to $\omega_{n+1}^4 = 1$. A simple calculation then gives us that $\omega_{n+1}^4 = 1$ for $k \in \C^+$ if and only if $k = i a_{n+1}$. After a lengthy computation using \eqref{Rn0}, we obtain that $R_n^0(k) = 0$ for $k \in \C^+$ if and only if equation \eqref{detM} holds.
 
Now suppose that $q_N(k) = 0$, $\det(M_{N-1}...M_2M_1)(k) \neq 0$ at $k = i\kappa$, and that $k$ is not a pole of $\widetilde{R_N^0}$. Then $k$ is a pole of $B_N$, and consequently a pole of $R = R_{1,...,N}$ as well. Consequently, $k^2$ is a bound state of the Schr\"odinger equation. Since $\det(M_{N-1}...M_2M_1)(k) \neq 0$, we have that $p_N(k) \neq 0$. Since $k$ is not a pole of $\widetilde{R_N^0}$ and since $R_N^0$ and $\widetilde{R_N^0}$ have the same zeros with the same multiplicity, $-\frac{R_N^0}{\widetilde{R_N^0}}p_n \neq 0$. However, $q_N = 0$ and the poles of $B_N$ are simple, so
\begin{equation}
	\Res_{k = i\kappa} B_N = -\frac{R_N^0}{\widetilde{R_N^0}} \frac{p_N}{q_N'}. \label{resB}
\end{equation}
To find $q_N'$, we can differentiate \eqref{pqrec} to acquire
\begin{equation}
\begin{pmatrix} p_{n+1}' \\ q_{n+1}' \end{pmatrix} = M_n \begin{pmatrix} p_n' \\ q_n' \end{pmatrix} + M_n' \begin{pmatrix} p_n \\ q_n \end{pmatrix}. \label{p'q'rec}
\end{equation}

Therefore, for the poles of $R$ where $\det(M_N...M_2M_1) \neq 0$ and that are not poles of $\widetilde{R_N^0}$, the residues of $R$ can be recovered through \eqref{resB} and \eqref{resR2}.

\section{Numerical Simulations}
Tables \ref{tab:bound_block}, \ref{tab:bound_wave1}, and \ref{tab:bound_wave2} give a comparison of some algorithms for calculating bound states described below. The exact bound states in table \ref{tab:bound_block} were calculated using equations \eqref{bound_oneblock}. All calculations were performed using MATLAB.

There are two commonly used numerical methods for approximating the bound states:
\begin{itemize}
\item[(1)] Matrix methods - Estimate the Schr\"odinger operator $H = -\frac{d^2}{dx^2} + V(x)$ using a finite-dimensional matrix and find the eigenvalues of the matrix. In particular, \cite{T00} describes how this can be done using the Fourier basis. In tables \ref{tab:bound_block}, \ref{tab:bound_wave1}, and \ref{tab:bound_wave2}, a $512 \times 512$ matrix is used.
\item[(2)] Shooting Method - The Shooting Method involves recursively choosing values of $\lambda$ and ``shooting'' from both end points
to a point $c \in [a,b]$. Define the miss-distance function to be the Wronskian determinant
\begin{equation*}
	D(\lambda) = \begin{vmatrix} u_L(c,\lambda) & u_R(c,\lambda)
														\\  u_L'(c,\lambda)  & u_R'(c,\lambda) \end{vmatrix}.
\end{equation*}
where $u_L$ is the interpolated function from the left endpoint and $u_R$ is from the right endpoint. If $\lambda$ is an eigenvalue that satisfies the boundary value problem, then $D(\lambda) = 0$. For more details, see for example \cite{P93}.
\end{itemize}
There are of course other existing methods for approximating bound states; see for example \cites{FV87, K07, KL81, C05}. However, we will only focus on these two.

If one approximates the potential using finitely many blocks, then we can use the following algorithms for estimating bound states:
\begin{itemize}
\item[(3)] Use the recursive formulas \eqref{recR} and \eqref{recA} to find the bound states as zeros of $1/R_{1,...,N}$.
\item[(4)] Similarly, one can use the recursive formula \eqref{recB} to find the bound states as zeros of $1/B_N$.
\item[(5)] Using \eqref{mateq}, the bound states can be found as zeros of $q_N$. One must also check the values of $k$ listed in \eqref{detM} where $\det(M_{N-1}...M_1) = 0$.
\end{itemize}

Algorithm (1) seems to be the fastest of these algorithms, followed closely by (2). Moreover, algorithm (1) has great accuracy when the initial potential is smooth. However, for discontinuous potentials, the Gibbs phenomenon severely hinders the accuracy of the algorithm. Moreover, the domain chosen seems to effect algorithms (1) and (2) greatly. On the other hand, tables \ref{tab:bound_wave1} and \ref{tab:bound_wave2} demonstrate that algorithms (3)-(5) are more robust when choosing the domain, with a smaller domain being preferable for the amount of time. All of algorithms (2)-(5) rely on finding roots of some function, so inheritently all of these functions have all of the problems that root finders tend to have. For example, given a good initial approximation of a bound state, the root finder might diverge or converge to a different bound state. Furthermore, the bound states are known to cluster towards 0, which makes it increasingly difficult to accurately determine all of the bound states as the number of bound states increases. However, when the root finders do converge, algorithms (3)-(5) are extremely accurate. Algorithms (3)-(5) also seem to be much slower than algorithms (1) and (2), with (5) being the slowest. 

In summary, the commonly used algorithms (1) and (2) for calculating bound states are much faster than the other algorithms. Moreover, algorithm (1) tends to be extremely accurate, especially when the potential is smooth. However, although algorithms (3)-(5) are much slower, they also tend to be very accurate, especially with discontinuous potentials. Moreover, these algorithms seem to be more robust when choosing the domain of the potential.

Supposing the bound states have been calculated, tables \ref{tab:norm_block} and \ref{tab:norm_wave1} give a comparison of some of the various algorithms for computing (left) norming constants. First is the algorithm described in the present paper:
\begin{itemize}
\item [(i)] The potential is approximated using finitely many blocks, and the norming constants are calculated as residues via equations \eqref{resB} and \eqref{resR2}.
\end{itemize}

Next we have the obvious algorithm using the definition of the left norming constant:
\begin{itemize}
\item[(ii)] Suppose $V$ has compact support $[A,B]$. We have that $\phi(x,k) = \phi_l(x,k)/T(k)$ satisfies $\phi(x,k) = e^{ikx}$ for $x \geq B$. One can numerically integrate the Schr\"odinger equation from $B$ to $A$. Then $c_l^2 = \|\phi\|_2^{-1}$, which can be numerically integrated.
\end{itemize}

The authors were also presented the following algorithms by Paul Sacks: letting $a = 1/T$ and $b = -R/T$, then $R = -\frac{b}{a}$ and the transition matrix $\Lambda$ given in \eqref{transition_matrix} becomes
\begin{equation*}
	\Lambda = \begin{pmatrix} a & b \\ \widetilde{b} & \widetilde{a} \end{pmatrix}.
\end{equation*}
Moreover, $b$ is analytic everywhere in $\C^+$, and the simple poles of $T$ in $\C^+$ are simple zeros of $a$. Consequently, \eqref{ResR} gives us that
\begin{equation*}
	c_n^2 = i \frac{b(i\kappa_n)}{a'(i\kappa_n)}.
\end{equation*}
The derivative $a'$ with respect to $k$ can be approximated using the central difference
\begin{equation*}
	a'(k) \approx \frac{a(k+\eta/2) - a(k-\eta/2)}{\eta}.
\end{equation*}
The question then becomes how one evaluates $a(k)$ and $b(k)$. Here are two approaches:
\begin{itemize}
\item[(iii)] The potential is approximated using a finite number of blocks, and $a$ and $b$ are calculated using potential fragmention \eqref{potentialfragmentation}. The transition matrices are evaluated using equation \eqref{oneblock}.
\item[(iv)] Supposing the potential has compact support $[\alpha,\beta]$, the Schrodinger equation can be numerically integrated from $\alpha$ to $\beta$ with the initial conditions $\phi(\alpha,k) = e^{-ik\alpha}$, $\phi'(\alpha,k) = -ike^{-ik\alpha}$. Then $\phi(x,k) = \phi_r(x,k)/T(k)$, so for $x \geq \beta$
\begin{equation*}
	\phi(x,k) = a(k)e^{-ikx} - b(k)e^{ikx}.
\end{equation*}
Consequently, $a$ and $b$ can be retrieved from
\begin{equation*}
	\begin{pmatrix} a(k) \\ b(k) \end{pmatrix} = \frac{1}{2} \begin{pmatrix} e^{ik\beta} & \frac{ie^{ik\beta}}{k} \\ -e^{-ik\beta} & \frac{ie^{-ik\beta}}{k} \end{pmatrix} \begin{pmatrix} \phi(\alpha,k) \\ \phi'(\alpha,k) \end{pmatrix}.
\end{equation*}
\end{itemize}

Algorithms (ii) and (iv) seem to be the fastest of these four algorithms. However, they are also the least accurate since they require integrating the Schr\"odinger equation which is extremely sensitive to errors. Algorithms (iii) generally takes about half as long as algorithm (i). Algorithm (i) seems to be the most accurate for discontinuous potentials, while algorithm (iii) seems to be the most accurate for smooth potentials. Moreover, the accuracy of algorithms (i) and (iii) increases when the bound states are approximated using algorithms (3)-(5).

\begin{figure}
	\includegraphics[width=4in]{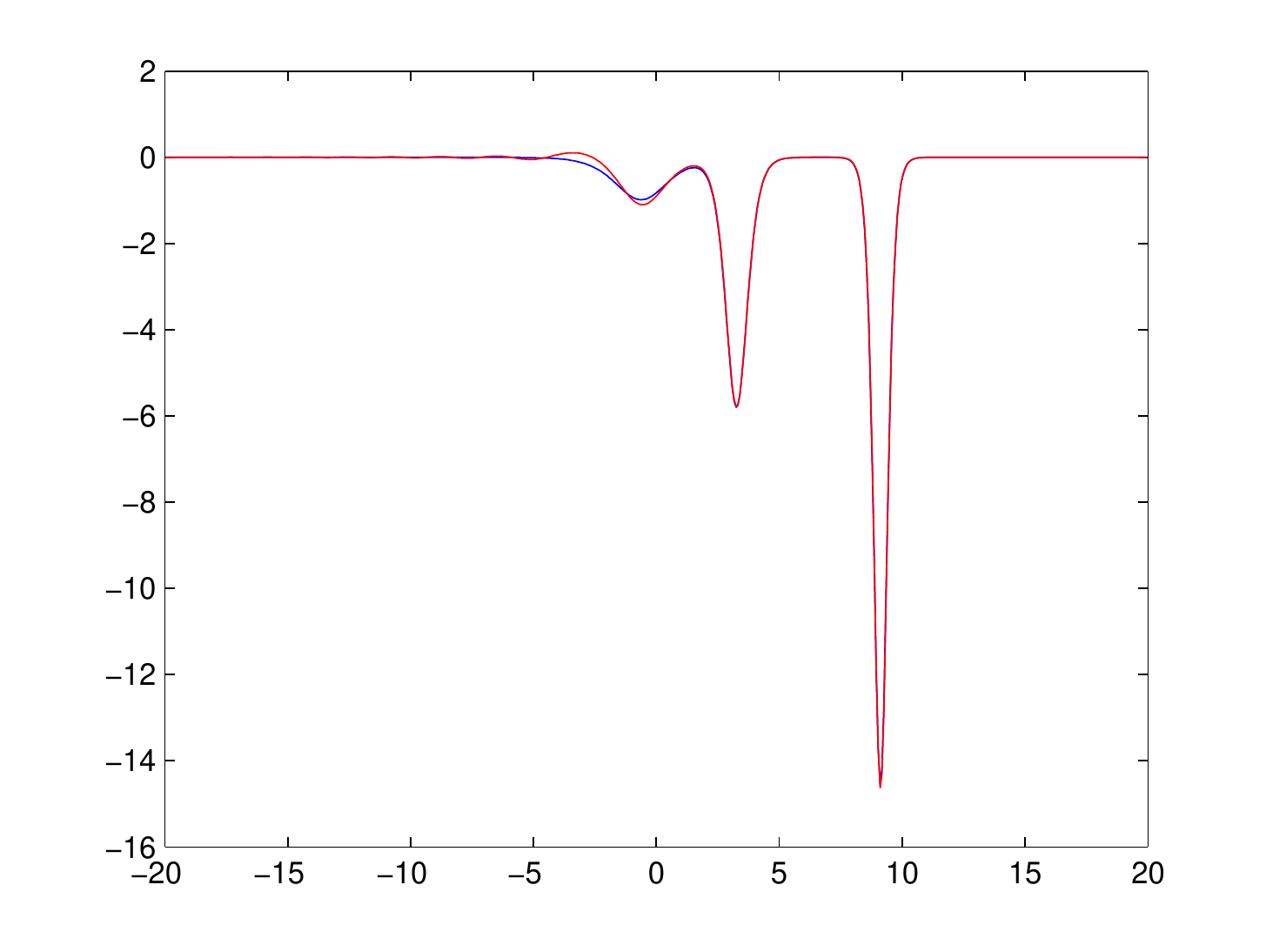}
	\caption{$V(x) = -10\sech^2(x)$, $t = 0.3$}	\label{fig:plot1}
\end{figure}

\begin{figure}
	\includegraphics[width=4in]{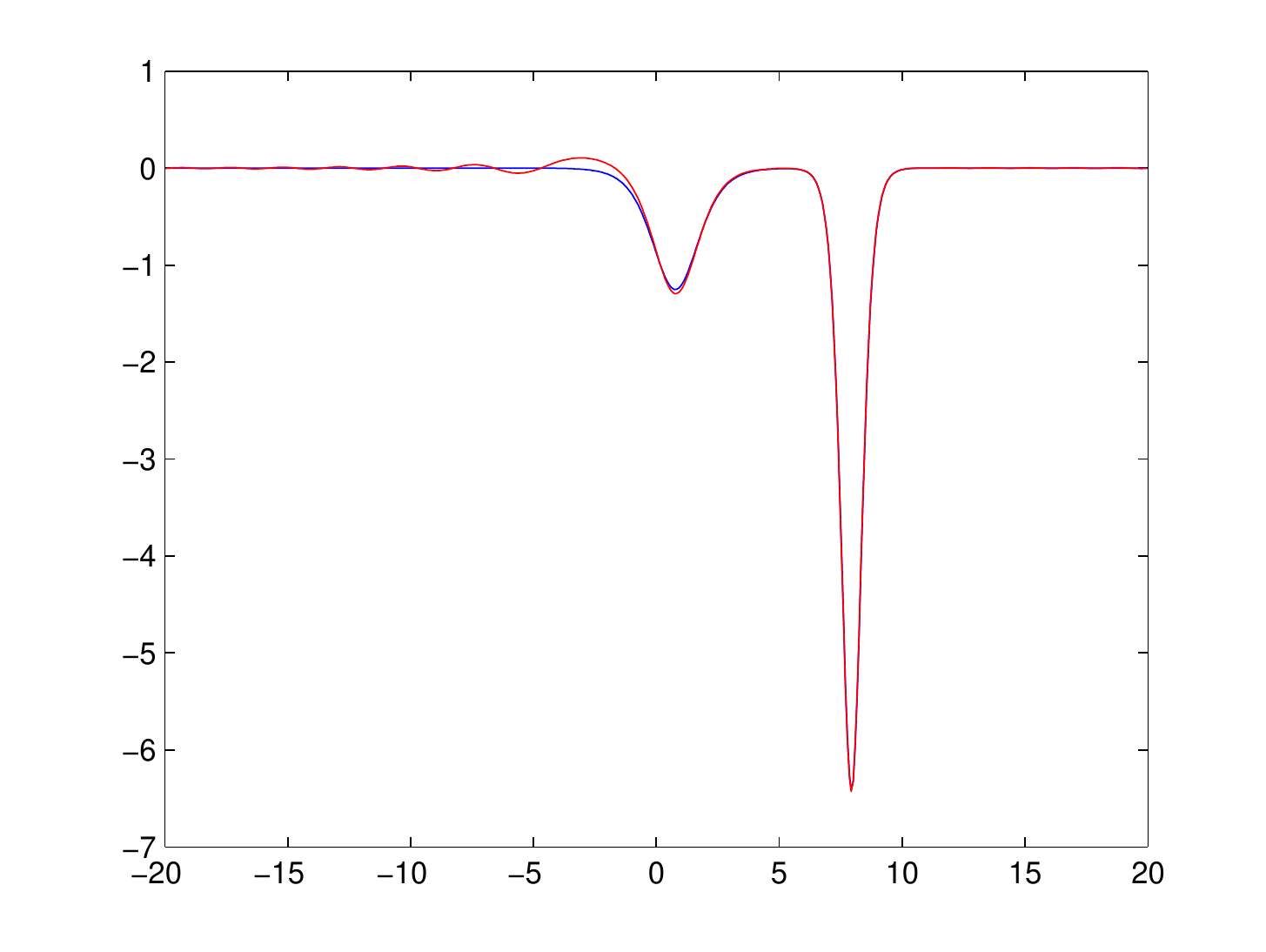}
	\caption{$V(x) = -5\sech^2(x)$, $t = 0.6$}	\label{fig:plot2}
\end{figure}

Lastly, figures \ref{fig:plot1} and \ref{fig:plot2} compare the asymptotic formula given in \cite{AC91} with the numerically integrated solution obtained by using the split step Fourier method.

\begin{table}\scriptsize	
\caption{$V(x) = -4\chi_{[-4,0]}(x)$, domain chosen $[-10,10]$, spacial step size $h = 0.01$} \label{tab:bound_block}
\begin{tabular}{rrrrrr}
Algorithm & $\kappa_1$      &  $\kappa_2$        &  $\kappa_3$        &  Relative Error	 &  Time (sec)	\\
\hline
Exact & 	 1.899448036751944 &  1.571342556813314 &  0.876610362727433 &                 0  &  0.004355000000000	\\
(1)	&   1.898826427139628 &  1.568514453040000 &  0.867505110670815 &  0.003651829842877 &  0.126239000000000	\\
(2)	&   1.899418261950639 &  1.572105829640451 &  0.872097420881459 &  0.001749410414267 &  0.505034000000000	\\
(3)	&   1.899448036751942 &  1.571342556813313 &  0.876610362727439 &  0.000000000000003 &  4.168762000000000	\\
(4)	&   1.899448036751949 &  1.571342556813312 &  0.876610362727428 &  0.000000000000003 &  5.425778000000000	\\
(5)	&   1.899448036751942 &  1.571342556813315 &  0.876610362727434 &  0.000000000000001 & 10.268152000000001	
\end{tabular}
\end{table}

\begin{table}\scriptsize	
\caption{$V(x) = -4\chi_{[-4,0]}(x)$, domain chosen $[-4,0]$, spacial step size $h = 0.01$, energy step size $\eta = 0.001$, exact bound states used} \label{tab:norm_block}
\begin{tabular}{rrrrrr}
Algorithm & $c_1^2$      	&  $c_2^2$        	&  $c_3^2$        	&  Relative Error	 	&  Time (sec)	\\
\hline
Exact &  0.038798932148319 &  0.145167980693995 &  0.257227284424067 &                  0 &  0.005992000000000	\\
(i)   &	0.038798932148326 &  0.145167980694058 &  0.257227284424741 &  0.000000000002272 &  2.008827000000000	\\
(ii)	&	0.141300713908832 &  0.293968570328614 &  0.444084025980906	&  0.872538777834092	&	0.032151000000000	\\
(iii) & 	0.038798937542783 &  0.145168027811526 &  0.257226712349713 &  0.000001926938416 &  2.070128000000000	\\
(iv) 	&  0.051311576782601 &  0.109225786002665 & -0.041977058580690 &  1.012467614318974 &  0.147137000000000
\end{tabular}
\end{table}

\begin{table}\scriptsize	
\caption{$V(x) = -\sech^2(x)$, domain chosen $[-5,5]$, spacial step size $h = 0.01$} \label{tab:bound_wave1}
\begin{tabular}{rrrr}
Algorithm & $\kappa$       &  Relative Error	 &  Time (sec)	\\
\hline
Exact &  1.000000000000000 &                  0 &                  0	\\
(1)   & 	1.000181385743159 &  0.000181385743159 &  0.123699000000000	\\
(2)   & 	1.000010661550817 &  0.000010661550817 &  0.165820000000000	\\
(3)   & 	0.999997769556372 &  0.000002230443628 &  8.624536000000001	\\
(4)   & 	0.999997769556371 &  0.000002230443629 &  8.780760000000001	\\
(5)   & 	0.999997769556372 &  0.000002230443628 & 14.264264000000001	
\end{tabular}
\end{table}

\begin{table}\scriptsize	
\caption{$V(x) = -\sech^2(x)$, domain chosen $[-5,5]$, spacial step size $h = 0.01$, energy step size $\eta = 0.001$, exact bound state used} \label{tab:norm_wave1}
\begin{tabular}{rrrr}
Algorithm & $c^2$       	&  Relative Error	 	&  Time (sec)	\\
\hline
Exact &  2.000000000000000 &                  0 &                  0	\\
(i)   &  2.004086813877857 &  0.002043406938928 &  3.460512000000000	\\
(ii)	&	1.274509474591987e-004 & 0.999936274526270 & 0.023956000000000	\\
(iii) &	1.999683571279579 &  0.000158214360211 &  1.631856000000000	\\
(iv) 	&  1.993946894122799 &  0.003026552938601 &  0.088449000000000	
\end{tabular}
\end{table}

\begin{table}\scriptsize	
\caption{$V(x) = -\sech^2(x)$, domain chosen $[-10,10]$, spacial step size $h = 0.01$} \label{tab:bound_wave2}
\begin{tabular}{rrrr}
Algorithm & $\kappa$       &  Relative Error	 &  Time (sec)	\\
\hline
Exact &  1.000000000000000 &                  0 &                  0	\\
(1)   & 	1.000000008244449 &  0.000000008244449 &  0.119652000000000	\\
(2)   & 	1.000071226867798 &  0.000071226867798 &  0.177487000000000	\\
(3)   & 	0.999997777799307 &  0.000002222200693 & 17.526257000000001	\\
(4)   & 	0.999997777799307 &  0.000002222200693 & 18.236691000000000	\\
(5)   & 	0.999997777799307 &  0.000002222200693 & 28.798037999999998	
\end{tabular}
\end{table}

\section{Haar Systems and a KdV Large-Time Solver } 
Suppose now that $V$ is finite, nonpositive, and has compact support. Then $V$ can be well approximated using finitely many nonpositive blocks. For such potentials $V$, we now summarize the algorithm for solving the KdV for large times:
\begin{itemize}
\item
Approximate the potential $V(x)$ using $N$ nonpositive blocks
\item
Bound states are found as zeros of $1/R_{1,...,N}$ with initial estimates, for example, derived from a spectral matrix estimate of the Sch\"odinger operator
\item
The norming constants are calculated as residues of $B_N$ at the bound states using the previously described recursive formulas
\item
The solution to the KdV is obtained from the asymptotic formula:
\begin{equation*}
	u(x,t) \sim -2 \sum_{n=1}^N \kappa_n^2 \sech^2(\kappa_nx - 4\kappa_n^3t + \ln\sqrt{\gamma_n})
\end{equation*}
where
\begin{equation*}
	\gamma_n = \frac{2\kappa_n}{c_n^2}\prod_{m=1}^{n-1}\left(\frac{\kappa_n+\kappa_m}{\kappa_n-\kappa_m}\right)^2.
\end{equation*}
\end{itemize}

There are a number of possible improvements to this algorithm. For example, the number of bound states is known for a single block, so the results in \cite{AKM98} could possibly be implemented to determine the exact number of bound states for the potential. As another example, instead of piecewise-constant functions, one could instead use higher order spline interpolants of the potential. All of the recursive formulas in section \ref{sec:multiblocks} were derived from potential fragmentation, which holds for arbitrary potentials; the only things that would change would be the formula for $R_n^0$, the initial values in the recursive formulas, and the values for $k$ in \eqref{detM}. For example, in the case of piecewise-linear spline interpolants, the formula for $R_n^0$ would involve the Airy functions.

Another possible route for improvement would be the use of Haar wavelets or other wavelets. We will only consider Haar wavelets in the current paper. For a great exposition on Haar and other wavelets, see \cite{O06}. Consider the \it{scaling function}
\begin{equation*}
	\varphi(x) = \varphi_0(x) = \begin{cases} 1 & \text{if } 0 < x \leq 1,	\\
													0 & \text{otherwise},
				 				 \end{cases}
\end{equation*}
and the \it{mother wavelet}
\begin{equation*}
	w(x) = \begin{cases} 1 & \text{if } 0 < x \leq 1/2,	\\
							  -1 & \text{if } 1/2 < x \leq 1,	\\
							  	0 & \text{otherwise}.
			 \end{cases}
\end{equation*}
We form the Haar wavelets as follows: let
\begin{equation*}
	w_{j,0}(x) = w(2^j x).
\end{equation*}
Then $w_{j,0}$ has support $[0,2^{-j}]$. Next, we translate $w_{j,0}$ so as to fill up the entire interval $[0,1]$ with $2^j$ subintervals of length $2^{-j}$:
\begin{equation*}
	w_{j,k}(x) = \varphi_{2^j+k} = w_{j,0}(x-k) = w(2^j(x-k)),	\	\	\	k = 0,1,...,2^j-1.
\end{equation*}
Then $w_{j,k}$ has support $[2^{-j}k,2^{-j}(k+1)]$. The collection of \it{Haar wavelets}
\begin{equation*}
	\mathcal{H}_{2^n} = \{ \varphi_m \ : \ 0 \leq m \leq 2^n-1 \}
\end{equation*}
forms an orthogonal system with respect to the $L^2$ norm of dimension $2^n$; the collection $\mathcal{H}_\infty$ forms a complete orthogonal system for $L^2([0,1])$. For $\mathcal{H}_{2^n}$, let $\mathbf{\varphi}_r$ denote the vector in $\R^{2^n}$ corresponding to $\varphi_r$; i.e., the entries of $\mathbf{\varphi}_r$ are the function values of $\varphi_r$ on the $2^n$ intervals.

By translating and scaling, suppose without loss of generality that $V$ has compact support $[0,1]$. Since $V$ is finite, we have that $V \in L^2([0,1])$, so $V$ can be expressed in terms of the Haar basis:
\begin{equation*}
		V = \sum_{r=0}^\infty c_r \varphi_r
\end{equation*}
where
\begin{equation*}
	c_r = \frac{\left\langle V, \varphi_r \right\rangle_2}{\|\varphi_r\|_2}.
\end{equation*}

Let $V_0$ denote the piecewise-constant approximation of $V$ on the $2^n$ intervals mentioned above, and let $\mathbf{V}$ denote the corresponding column vector in $\R^{2^n}$. Then $V_0$ can be represented as a linear combination of the Haar wavelets in $\mathcal{H}_{2^n}$:
\begin{equation*}
	V_0 = \sum_{r=0}^{2^n-1} c_r \varphi_r
\end{equation*}
where the coefficients $c_r$ are as described above. Letting $\mathbf{c}$ denote the column vector of coefficients $c_r$, the \it{discrete wavelet transform} (DWT) is the map $H_{2^n} : \mathbf{V} \mapsto \mathbf{c}$; that is, $H_{2^n}$ is a change of basis from the standard basis to the Haar basis. Letting $W_{2^n}$ denote the matrix whose $r$-th column is $\mathbf{\varphi}_r$, we have that
\begin{equation*}
	\mathbf{V} = W_{2^n} \mathbf{c},
\end{equation*}
so
\begin{equation*}
	\mathbf{c} = W_{2^n}^{-1} \mathbf{V},
\end{equation*}
implying that $H_{2^n} = W_{2^n}^{-1}$. (Note: often, the columns are normalized so that $W_{2^n}$ is an orthogonal matrix. In this case, $H_{2^n} = W_{2^n}^*$ where $*$ denotes the transpose).

The Discrete Wavelet Transform is analogous to the Fast Fourier Transform (FFT), which expresses $\mathbf{V}$ in the orthogonal basis corresponding to the Fourier basis $\{ e^{i2^n x} : -2^{n-1} < r \leq 2^n \}$ in $L^2([-\pi,\pi])$. However, the Fourier basis is not localized, unlike the Haar basis, so the Fourier basis has difficulty capturing data concentrated in a relatively small region. The Fourier basis tends to accurately approximate smoother functions, while exhibiting the so called \it{Gibb's phenomenon} at discontinuities. On the other hand, the Haar basis tends to accurately approximate discontinuous functions, while only slowly converging to smoother functions.

In the context of solving the KdV, Haar wavelets may possibly be implemented in a couple ways. One approach would be to approximate the potential using Haar wavelets since it generally gives more accurate piecewise-constant interpolants than, say the midpoint rule. Then the interpolating potential would be changed to the standard basis and used in our algorithm. 

Currently, our potentials are being approximated by step functions using the standard basis since this is the form required for potential fragmentation. However, it is more desirable to represent the potential using Haar wavelets in many cases, such as for signal processing. Another approach for improving the algorithm would be to change all of our formulas over to the Haar basis. There are still many open questions in this regard, for example
\begin{itemize}
\item[(I)] If our potential was to be expressed in the Haar basis instead of the standard basis, what would be an efficient way to determine the scattering data? 
\item[(II)] Could potential fragmentation and our recursive formulas be modified to use the Haar representation of the potential?
\end{itemize}

%\section{Conclusions}

\section{Acknowledgements}
This work was done as part of the REU program run by the third author in the summer of 2009, and was supported by NSF grants DMS 0707476 and DMS 1009673. The government support is highly appreciated.

We would like to thank the other participants who have also contributed to this project: Lyman Gillispie and Sigourney Walker. We would particularly like to thank Paul Sacks for providing us algorithms (iii) and (iv) for calculating norming constants. We are also grateful to Constantine Khroulev and Anton Kulchitsky for useful consultations and discussions.

\end{document}